\documentclass[12pt, reqno]{amsart}
\usepackage{amsmath, amsthm, amscd, amsfonts, amssymb, graphicx, color}
\usepackage[bookmarksnumbered, colorlinks, plainpages]{hyperref}

\newtheorem{theorem}{Theorem}[section]
\newtheorem{lemma}[theorem]{Lemma}

\theoremstyle{definition}
\newtheorem{definition}[theorem]{Definition}

\theoremstyle{remark}

\numberwithin{equation}{section}

\begin{document}
\setcounter{page}{1}

\title[Reversible AJW-algebras]{Reversible AJW-algebras}

\author[Sh. Ayupov]{Shavkat Ayupov$^1$}

\address{$^{1}$ Institute of Mathematics, National University of Uzbekistan,
Tashkent, Uzbekistan.} \email{\textcolor[rgb]{0.00,0.00,0.84}{
sh$_-$ayupov@mail.ru}}

\author[F. Arzikulov]{Farhodjon Arzikulov$^2$}

\address{$^{2}$ Department of Mathematics, Andizhan State University, Andizhan, Uzbekistan.}
\email{\textcolor[rgb]{0.00,0.00,0.84}{arzikulovfn@rambler.ru}}




\date{Received: xxxxxx 14.09.26; Revised: yyyyyy; Accepted: zzzzzz.
\newline \indent Partially supported by TWAS, The Abdus Salam, International
Centre, for Theoretical Physics (ICTP), Grant:
13-244RG/MATHS/AS$_-$I-UNESCO FR:3240277696}

\begin{abstract}
In this article it is proved that for every special AJW-algebra
$A$ there exist central projections $e$, $f$, $g\in A$, $e+f+g=1$
such that (1) $eA$ is reversible and there exists a norm-closed
two sided ideal $I$ of $C^*(eA)$ such that
$eA={{}^\perp}(^\perp(I_{sa})_+)_+$; (2) $fA$ is reversible and
$R^*(fA)\cap iR^*(fA)=\{0\}$; (3) $gA$ is a totally nonreversible
AJW-algebra.
\end{abstract}

\maketitle

\section{Introduction}

This article is devoted to abstract Jordan operator algebras,
which are analogues of abstract W$^*$-algebras (AW$^*$-algebras)
of Kaplansky. These Jordan operator algebras can be characterized
as a JB-algebra satisfying the following conditions

(1)~in the partially ordered set of all projections any subset of
pairwise orthogonal projections has the least upper bound in this
JB-algebra;

(2)~every maximal associative subalgebra of this JB-algebra is
generated by it's projections $($i.e. coincides with the least
closed subalgebra containing all projections of the given
subalgebra $)$.

In the articles \cite{Arz1}, \cite{Arz2} the second author
introduced analogues of annihilators for Jordan algebras and gave
algebraic conditions equivalent to (1) and (2). Currently, these
JB-algebras are called AJW-algebras or Baer JB-algebras in the
literature. Further, in \cite{Arz3} a classification of these
algebras has been obtained. It should be noted that many of facts
of the theory of JBW-algebras and their proofs hold for
AJW-algebras. For example, similar to a JBW-algebra an AJW-algebra
is the direct sum of special and purely exceptional Jordan
algebras \cite{Arz3}.

It is known from the theory of JBW-algebras that every special
JBW-algebra can be decomposed into the direct sum of totally
irreversible and reversible subalgebras. In turn, every reversible
special JBW-algebra decomposes into a direct sum of subalgebra,
which is the hermitian part of a von Neumann algebra and the
subalgebra, enveloping real von Neumann algebra of which is purely
real \cite{ASh}, \cite{Ay2}. In this paper we prove a similar
result for AJW-algebras, the proof of which requires a different
approach. Namely, we prove that for every special AJW-algebra $A$
there exist central projections $e$, $f$, $g\in A$, $e+f+g=1$ such
that (1) $eA$ is reversible and there exists a norm-closed two
sided ideal $I$ of $C^*(eA)$ such that
$eA={{}^\perp}(^\perp(I_{sa})_+)_+$; (2) $fA$ is reversible and
$R^*(fA)\cap iR^*(fA)=\{0\}$; (3) $gA$ is a totally nonreversible
AJW-algebra.

\section{Preliminary Notes}

We fix the following terminology and notations.

Let $\mathcal{A}$ be a real Banach $*$-algebra. $\mathcal{A}$ is
called a real C$^*$-algebra, if
$\mathcal{A}_c=\mathcal{A}+i\mathcal{A}=\{a+ib: a,b\in
\mathcal{A}\}$, can be normed to become a (complex) C$^*$-algebra,
and keeps the original norm on $\mathcal{A}$ \cite{3}.

Let $A$ be a JB-algebra, $P(A)$ be a set of all projections of
$A$. Further we will use the following standard notations:
$\{aba\}=U_a b:=2a(ab)-a^2b$, $\{abc\}=a(bc)+(ac)b-(ab)c$ and
$\{aAb\}=\{\{acb\}:c\in A\}$, where $a$, $b$, $c\in A$. A
JB-algebra $A$ is called an AJW-algebra, if the following
conditions hold:

(1)~in the partially ordered set $P(A)$ of projections any subset
of pairwise orthogonal projections has the least upper bound in
$A$;

(2)~every maximal associative subalgebra $A_o$ of the algebra $A$
is generated by it's projections $($i.e. coincides with the least
closed subalgebra containing $A_o\cap P(A))$.

Let $(S)^\perp=\{a\in A:\,(\forall\, x\in S)\ U_{a}x=0\}$,\ \
${}^{\perp}(S)=\{x\in A:\,(\forall\, a\in S)\ U_{a}x=0\}$,\ \
$^{\perp}(S)_+={^{\perp}}(S)\cap A_+$.

Then for a JB-algebra $A$ the following conditions are equivalent:

(1) $A$ is an AJW-algebra;

(2) for every subset $S\subset A_+$ there exists a projection
$e\in A$ such that $(S)^\perp =U_{e}(A)$;

(3) for every subset $S\subset A$ there exists a projection $e\in
A$ such that $^{\perp }(S)_+=U_{e}(A_+)$ \cite{Arz1}.

Let $A$ be a real or complex $*$-algebra, and let $S$ be a
nonempty subset of $A$. Then the set $R(S)=\{x\in A: sx=0\, for
\,\,all\, s\in S\}$ is called the right annihilator of $S$ and the
set $L(S)=\{x\in A: xs=0\,\, for \,\,all\,\,s\in S\}$ is called
the left annihilator of $S$. A $*$-algebra $A$ is called a Baer
$*$-algebra, if the right annihilator of any nonempty set
$S\subseteq A$ is generated by a projection, i.e. $R(S)=gA$ for
some projection $g\in A$ ($g^2=g=g^*$). If $S=\{a\}$ then the
projection $1-g$ such that $R(S)=gA$ is called the right
projection and denoted by $r(a)$. Similarly one can define the
left projection $l(a)$. A (real) C$^*$-algebra $A$, which is a
Baer (real) $*$-algebra, is called an (real) AW$^*$-algebra
\cite{5}, \cite{6}. Real AW$^*$-algebras were introduced and
investigated in \cite{6}, \cite{10}. In these papers it was shown
that for a real AW$^*$-algebra $\mathcal{A}$ the C$^*$-algebra
$M=\mathcal{A}+i\mathcal{A}$ is not necessarily a complex
AW$^*$-algebra.

Let $A$ be an AJW-algebra. By \cite[Theorem 2.3]{Arz3} we have the
equality $A=A_I\oplus A_{II}\oplus A_{III}$, where $A_I$ is an
AJW-algebra of type I, $A_{II}$ is an AJW-algebra of type II and
$A_{III}$ is an AJW-algebra of type III \cite{Arz3}. By
\cite[Theorem 3.7]{Arz3} $A_I$, in its turn, is a direct sum of
the following form
$$
A_I=A_\infty \oplus A_1\oplus A_2\oplus \dots,
$$
where $A_n$ for every $n$ either is $\{0\}$ or an AJW-algebra of
type I$_n$, $A_\infty$ is a direct sum of AJW-algebras of type
I$_\alpha$ with $\alpha$ infinite. If $A=A_1\oplus A_2\oplus
\dots$ then $A$ is called an AJW-algebra of type I$_{fin}$ and
denoted by $A_{I_{fin}}$ and if $A=A_\infty$ then $A$ is called an
AJW-algebra of type I$_\infty$ and denoted by $A_{I_\infty}$. We
say that $A$ is properly infinite if $A$ has no nonzero central
modular projection. The fact that an AJW-algebra $A_{II}$ of type
II is a JC-algebra can be proved similar to JBW-algebras
\cite{HOS}. Therefore, it is isomorphic to some AJW-algebra
defined in \cite{Top} (i.e. to some AJW-algebra of self-adjoint
operators), and by virtue of \cite{Top} $A_{II}=A_{II_1}\oplus
A_{II_\infty}$, where $A_{II_1}$ is a modular AJW-algebra of type
II and $A_{II_\infty}$ is an AJW-algebra of type II, which is
properly infinite. So, we have the decomposition
$$
A=A_{I_{fin}}\oplus A_{I_\infty}\oplus A_{II_1}\oplus
A_{II_\infty}\oplus A_{III}.
$$
It is easy to verify that the part  $A_{I_{fin}}\oplus A_{II_1}$
is modular, and $A_{I_\infty}\oplus A_{II_\infty}\oplus A_{III}$
is properly infinite (i.e. properly nonmodular).

\section{Reversibility of AJW-algebras}

Let $A$ be a special AJW-algebra on a complex Hilbert space $H$.
By $R^*(A)$ we denote the uniformly closed real $*$-algebra in
$B(H)$, generated by $A$, and by $C^*(A)$ the C$^*$-algebra,
generated by $A$. Thus the set of elements of kind
$$
\sum_{i=1}^n  \prod_{j=1}^{m_i} a_{ij} (a_{ij}\in A)
$$
is uniformly dense in $R^*(A)$. Let $iR^*(A)$ be the set of
elements of kind $ia$, $a\in R^*(A)$. Then
$C^*(A)=R^*(A)+iR^*(A)$. \cite{ASh}, \cite{SE}

\begin{lemma}   \label{2.1}
The set $R^*(A)\cap iR^*(A)$ is a uniformly closed two sided ideal
in $C^*(A)$.
\end{lemma}

{\it Proof.} If $a$, $b\in R^*(A)$ and $c=id\in R^*(A)\cap
iR^*(A)$, then $(a+ib)c=ac+ibid=ac-bd\in R^*(A)$. Similarly
$(a+ib)c\in iR^*(A)$, i.e. $(a+ib)c\in R^*(A)\cap iR^*(A)$. Since
$R^*(A)\cap iR^*(A)$ is uniformly closed and the set of elements
of kind $a+ib$, $a$, $b\in A$ is uniformly dense in $C^*(A)$, we
have $R^*(A)\cap iR^*(A)$ is a left ideal in $C^*(A)$. By the
symmetry $R^*(A)\cap iR^*(A)$ is a right ideal. $\triangleright$

Let $R$ be a $*$-algebra, $R_{sa}$ be the set of all self-adjoint
elements of $R$, i.e. $R_{sa}=\{a\in R: a^*=a\}$.

\begin{definition} A JC-algebra $A$ is said to be reversible
if $a_1a_2\dots a_n+a_na_{n-1}\dots a_1\in A$ for all $a_1,
a_2,\dots, a_n\in A$.
\end{definition}

Similar to JW-algebras we have the following criterion.

\begin{lemma} \label{2.3}
An AJW-algebra $A$ is reversible if and only if $A=R^*(A)_{sa}$.
\end{lemma}

{\it Proof.} It is clear that, if $A=R^*(A)_{sa}$, then $A$ is
reversible since
$$
(\prod_{i=1}^na_i+\prod_{i=n}^1a_i)^*=\prod_{i=n}^1a_i+\prod_{i=1}^na_i\in
R^*(A)_{sa}=A,
$$
for all $a_1$, $a_2$, $\dots$, $a_n\in A$. Conversely, let $A$ be
a reversible AJW-algebra. The inclusion $A\subset R^*(A)_{sa}$ is
evident. If $a=\sum_{i=1}^n\prod_{j=1}^{m_i}a_{ij}\in
R^*(A)_{sa}$, then
$$
a=\frac{1}{2}(a+a^*)=\frac{1}{2}\sum_{i=1}^n(\prod_{j=1}^{m_i}a_{ij}+\prod_{j=m_i}^1a_{ij})\in
A.
$$
Hence the converse inclusion holds, i.e. $R^*(A)_{sa}=A$.
$\triangleright$

\begin{lemma} \label{2.4}
Let $A$ be an AJW-algebra and let $I$  be a norm-closed ideal of
$A$. Then there exists a central projection $g$ such that
$^\perp(^\perp(I_{sa})_+)_+=gA_+$.
\end{lemma}

{\it Proof.} Since $A$ is an AJW-algebra there exists a projection
$g$ in $A$ such that
$$
^\perp(I_{sa})_+=U_{(1-g)}(A_+),
{{}^\perp}(^\perp(I_{sa})_+)_+=U_g(A_+),
$$
where $^\perp(S)_+=\{x\in A_+:(\forall a\in S) U_ax=0\}$ for
$S\subseteq A$.

Let $(u_\lambda)$ be an approximate identity of the JB-subalgebra
$I$ and $a$ be an arbitrary positive element in $I$. Then there
exists a maximal associative subalgebra $A_o$ of $A$ containing
$a$. Let $v_\mu$ be an approximate identity of $A_o$. Then
$(v_\mu)\subseteq (u_\lambda)$ and $\Vert av_\mu-a\Vert\to 0$. Let
$b\in A_+$ and
$$
U_{v_\mu}b=0
$$
for every $\mu$. Then $U_aU_{v_\mu}b=U_{av_\mu}b=0$ and
$U_cU_{av_\mu}b=0$, where $c$ is an element in $A$ such that
$b=c^2$. Hence $U_cU_{av_\mu}c^2=0$, $(U_c(av_\mu))^2=0$,
$U_c(av_\mu)=0$ and $U_cU_c(av_\mu)=U_b(av_\mu)=0$ for every
$\mu$. We have
$$
\Vert U_b(av_\mu)-U_ba\Vert=\Vert U_b(av_\mu)-a)\Vert\to 0
$$
because $\Vert av_\mu-a\Vert\to 0$ and the operator $U_b$ is
norm-continuous. Hence $U_ba=0$. We may assume that $a=d^2$ for
some element $d\in A$. Then
$$
U_dU_ba=U_dU_bd^2=(U_db)^2=0, U_db=0.
$$
Thus $U_dU_db=U_{d^2}b=U_ab=0$. Therefore, if $b\in
^\perp((u_\lambda))_+$ then $b\in ^\perp(I_{sa})_+$. Hence
$^\perp((u_\lambda))_+\subseteq  ^\perp(I_{sa})_+$. It is clear
that $^\perp(I_{sa})_+\subseteq ^\perp((u_\lambda))_+$ and
$$
^\perp(I_{sa})_+=^\perp((u_\lambda))_+.
$$
This implies that $^\perp((u_\lambda))_+=U_{(1-g)}(A_+)$ and
$$
\sup_\lambda u_\lambda=g.
$$

Let us prove that $U_g(A)$ is an ideal of $A$. Indeed, let $x$ be
an arbitrary element in $A$. Then $U_xu_\lambda\in I_{sa}$, i.e.
$U_xu_\lambda\in U_g(A)$. By \cite[Proposition 3.3.6]{HOS} and the
proof of \cite[Lemma 4.1.5]{HOS} we have $U_x$ is a normal
operator in $A$. Hence
$$
\sup_\lambda U_xu_\lambda=U_x(\sup_\lambda u_\lambda)=U_xg.
$$
At the same time
$$
\sup_\lambda U_xu_\lambda \in U_g(A).
$$
Hence $U_xg\in U_g(A)$. By \cite[2.8.10]{HOS} we have
$$
4(xg)^2=2gU_xg+U_xg^2+U_gx^2=2gU_xg+U_xg+U_gx^2.
$$
Therefore $(xg)^2\in U_g(A)$ and $xg\in U_g(A)$.

Now, let $y$ be an arbitrary element in $U_gA$. Then $y=U_gy$ and
$$
xy=(U_gx+\{gx(1-g)\}+U_{1-g}x)U_gy=U_gxU_gy+\{gx(1-g)\}U_gy\in
U_gA
$$
since $\{gx(1-g)\}\in U_gA$. Hence $U_gA$ is a norm-closed ideal
of $A$. Therefore $\{gA(1-g)\}=\{0\}$ and
$$
A=U_gA\oplus U_{1-g}A.
$$
This implies that $g$ is a central projection in $A$ and
${{}^\perp}(^\perp(I_{sa})_+)_+=gA_+$. $\triangleright$

\begin{lemma} \label{2.5}
Let $A$ be a reversible AJW-algebra on a Hilbert space $H$. Then
there exist two central projections $e$, $f$ in $A$ and a
norm-closed two sided ideal $I$ of $C^*(A)$ such that $e+f=1$,
$eA={{}^\perp}(^\perp(I_{sa})_+)_+$ and $R^*(fA)\cap
iR^*(fA)=\{0\}$.
\end{lemma}

{\it Proof.} Let $I=R^*(A)\cap iR^*(A)$. Since $A$ is reversible
by proposition \ref{2.3} we have $I_{sa}\subseteq A$. By
\cite[3.1]{ASh} $I$ is a two sided ideal of $C^*(A)$. Hence
$I_{sa}$ is an ideal of the AJW-algebra $A$. By proposition
\ref{2.4} we have there exists a central projection $g$ such that
${{}^\perp}(^\perp(I_{sa})_+)_+=gA_+$. It is clear that $g$ is a
central projection also in $C^*(A)$.

By the definitions of $I$ and $g$ we have
$$
R^*((1-g)A)\cap iR^*((1-g)A)=\{0\}.
$$
$\triangleright$

\begin{lemma} \label{2.6}
Let $A$ be an AJW-algebra and let $J$ be the set of elements $a\in
A$ such that $bac+c^*ab^*\in A$ for all $b$, $c\in C^*(A)$. Then
$J$ is a norm-closed ideal in $A$. Moreover $J$ is a reversible
AJW-algebra.
\end{lemma}

{\it Proof.} Let $a$, $b\in J$, $s$, $t\in C^*(A)$. Then
$$
s(a+b)t+t^*(a+b)s^*=(sat+t^*as^*)+(sbt+t^*bs^*)\in A,
$$
i.e. $J$ is a linear subspace of $A$. Now, if $a\in J$, $b\in A$,
$s$, $t\in C^*(A)$, then
$$
s(ab+ba)t+t^*(ab+ba)s^*=(sa(bt)+(bt)^*as^*)+
$$
$$
+((sb)at+t^*a(sb)^*)\in A,
$$
i.e. $J$ is a norm-closed ideal of $A$.

Let $a_1\in J$, $a_2$, $\dots$, $a_n\in A$ and $a=\prod_{i=2}^n
a_i$. Then $a_1a+a^*a_1\in A$ by the definition of $J$. Let us
show that $a_1a+a^*a_1\in J$; then, in particular, in the case of
$a_2$, $\dots$, $a_n\in J$ this will imply that $J$ is reversible.
For all $b$, $c\in C^*(A)$ we have
$$
b(a_1a+a^*a_1)c+c^*(a_1a+a^*a_1)b^*=
$$
$$
=(ba_1(ac)+(ac)^*a_1b^*)+((ba^*)a_1c+c^*a_1(ba^*)^*)\in A,
$$
i.e. $a_1a+a^*a_1\in J$. $\triangleright$

\begin{definition}
An AJW-algebra $A$ is said to be {\it totally nonreversible}, if
the ideal $J$ in lemma \ref{2.6} is equal to $\{0\}$, i.e.
$J=\{0\}$.
\end{definition}

\begin{theorem} \label{2.8}
Let $A$ be a special AJW-algebra. Then there exist central
projections $e$, $f$, $g\in A$, $e+f+g=1$ such that

(1) $J=(e+f)A$, $J$ is the ideal from lemma \ref{2.6};

(2) $eA$ is reversible and there exists a norm-closed two sided
ideal $I$ of $C^*(eA)$ such that
$eA={{}^\perp}(^\perp(I_{sa})_+)_+$;

(3) $fA$ is reversible and $R^*(fA)\cap iR^*(fA)=\{0\}$;

(4) $gA$ is a totally nonreversible AJW-algebra and
$$
gA=\sum_{\omega\in \Omega} C(Q_\omega,{\bf R}\oplus H_\omega),
$$
where $\Omega$ is a set of indices, $\{Q_\omega\}_{\omega\in
\Omega}$ is an appropriate family of extremal compacts and
$\{H_\omega\}_{\omega\in \Omega}$ is a family of Hilbert spaces.
\end{theorem}

{\it Proof.} We have
$$
A=A_1\oplus A_2\oplus \dots\oplus A_{I_\infty}\oplus
A_{II_1}\oplus A_{II_\infty}\oplus A_{III}
$$
and the subalgebra (without the part $A_2$)
$$
A_1\oplus A_3\oplus A_4\oplus\dots\oplus A_{I_\infty}\oplus
A_{II_1}\oplus A_{II_\infty}\oplus A_{III}
$$
is reversible. The last statement can be proven similar to
\cite[Theorem 5.3.10]{HOS}. By \cite[]{Kus} the subalgebra $A_2$
can be represented as follows
$$
A_2=\sum_{i\in \Xi} C(X_i,{\bf R}\oplus H_i),
$$
where $\Xi$ is a set of indices, $\{X_i\}_{i\in \Xi}$ is a family
of extremal compacts and $\{H_i\}_{i\in \Xi}$ is a family of
Hilbert spaces. Hence by \cite[Theorem 6.2.5]{HOS} there exist
central projections $h$, $g$ such that $A=hA\oplus gA$, $hA$ is
reversible and $gA$ is totally nonreversible. For all $a$,
$b_1$,$\dots$, $b_n$, $c_1$,$\dots$, $c_m$ in $hA$ we have
$$
b_1\dots b_n a c_1\dots c_m+c_mc_{m-1}\dots c_1 a b_nb_{n-1}\dots
b_1\in hA
$$
since $hA$ is reversible. Similarly for all $b$, $c$ in $R^*(hA)$,
$a\in hA$ we have
$$
bac+c^*ab^*\in hA.
$$
Hence $hA=J$

By proposition \ref{2.5} there exist two central projections $e$,
$f$ in $hA$ and a norm-closed two sided ideal $I$ of $C^*(hA)$
such that $e+f=h$, $eA={{}^\perp}(^\perp(I_{sa})_+)_+$, $fA$ is a
reversible AJW-algebra and $R^*(fA)\cap iR^*(fA)=\{0\}$. This
completes the proof. $\triangleright$

Let $A$ be a special AJW-algebra. Despite the fact that for the
real AW$^*$-algebra $R^*(A)$ the C$^*$-algebra $\mathcal{M}=
R^*(A)+iR^*(A)$ is not necessarily a complex AW$^*$-algebra we
consider, that

{\it Conjecture.} Under the conditions of theorem \ref{2.8} the
following equality is valid
$$
eA=I_{sa}.
$$

\end{document}